\documentclass[12pt]{article}

\usepackage{amsmath, epsfig, cite}
\usepackage{amssymb}
\usepackage{amsfonts}
\usepackage{latexsym}
\usepackage{amsthm}

\makeatletter
\renewcommand{\@seccntformat}[1]{{\csname the#1\endcsname}.\hspace{.5em}}
\makeatother

\newtheorem{thm}{Theorem}

\newtheorem{lem}[thm]{Lemma}

\newcommand{\pf}{\noindent{\it Proof.} }

\def\0{{\bf 0}}

\renewcommand{\qed}{\hfill$\Box$\medskip}

\begin{document}

\begin{center}
{\Large\bf A new proof of a theorem of Mansour and Sun}
\end{center}

\vskip 2mm \centerline{\large Victor J. W. Guo}

\begin{center}
{\footnotesize Department of Mathematics, East China Normal University,
Shanghai 200062, People's Republic of China\\
{\tt jwguo@math.ecnu.edu.cn\\
http://math.ecnu.edu.cn/\textasciitilde{jwguo}}}
\end{center}


\vskip 0.7cm \noindent{\small{\bf Abstract.} We give a new proof of
a theorem of Mansour and Sun by using number theory and Rothe's identity.


%


\vskip 0.7cm

It is well-known that the number of ways of choosing $k$ points, no
two consecutive, from a collection of $n$ points arranged on a cycle
is $\frac{n}{n-k}{n-k\choose k}$ (see \cite[Lemma 2.3.4]{Stanley}).
A generalization of this result was obtained by Kaplansky
\cite{Kaplansky}, who proved that the number of $k$-subsets
$\{x_1,\ldots,x_k\}$ of $\mathbb Z_n$ such that
$|x_i-x_j|\notin\{1,2,\ldots,p\}$ {\rm(}$1\leq i<j\leq k${\rm)} is
$\frac{n}{n-pk}{n-pk\choose k}$, where $n\geq pk+1$. Some other generalizations and related problems were
studied by several authors (see \cite{Chu,KP,Konvalina,Prodinger}). Very recently,
Mansour and Sun \cite{MS} extended Kaplansky's result as
follows.

\begin{thm}\label{thm:msun}
Let $m,p,k\geq 1$ and $n\geq mpk+1$. Then the number of $k$-subsets
$\{x_1,\ldots,x_k\}$ of $\mathbb Z_n$ such that
$|x_i-x_j|\notin\{m,2m,\ldots,pm\}$ for all $1\leq i<j\leq k$, denoted by $f_{m,n}$, is given by
$\frac{n}{n-pk}{n-pk\choose k}$.
\end{thm}

Their proof needs to establish a recurrence relation and compute the
residue of a Laurent series. Mansour and Sun \cite{MS} also asked for a combinatorial
proof of Theorem \ref{thm:msun}. In this note, we shall give a new but not purely combinatorial
proof of Theorem~\ref{thm:msun}.  Let $p$ and $k$ be fixed throughout.
Let $(a,b)$ denote the greatest common divisor of the
integers $a$ and $b$. We first establish the following three lemmas.

\begin{lem}\label{lem:one}
Let $(a,m)=1$ and let $d$ be a positive integer. Then at least one
of $a,\,a+m,\,a+2m,\ldots,a+(d-1)m$ is relatively prime to $d$.
\end{lem}
\pf If $(a,d)=1$, we are done. Now assume that
$(a,d)=p_1^{r_1}\ldots p_s^{r_s}$ and $d=p_1^{l_1}\cdots p_t^{l_t}$,
where $1\leq s\leq t$ and $p_1,\ldots,p_t$ are distinct primes and
$r_1,\ldots,r_s,l_1,\ldots,l_t\geq 1$. We claim that
$a+p_{s+1}\cdots p_t m$ is relatively prime to $d$. Indeed, since
$(a,m)=1$, we have $(p_1\cdots p_s,m)=1$ and therefore
\begin{align*}
(p_1\cdots p_s,a+p_{s+1}\cdots p_t m)    &=(p_1\cdots p_s,p_{s+1}\cdots p_t m)=1,\\
(p_{s+1}\cdots p_t,a+p_{s+1}\cdots p_t m)&=(p_{s+1}\cdots p_t,a)=1.
\end{align*}
This completes the proof. \qed

\begin{lem}\label{lem:two}
Let $(m,n)=d$. Then there exist integers $a,b$ such that $(a,n)=1$ and $am+bn=d$.
\end{lem}
\pf Since $(m,n)=d$, we may write $m=m_1d$ and $n=n_1d$, where $(m_1,n_1)=1$.
Then there exist integers $a$ and $b$ such that $am_1+bn_1=1$. It is clear that $(a,n_1)=1$.
Noticing that $(a+n_1t)m_1+(b-m_1t)n_1=1$, by Lemma \ref{lem:one},
we may assume that $(a,d)=1$ and so $(a,n)=1$.
\qed

\begin{lem}\label{lem:three}
Let $m,n\geq 1$ and $(m,n)=d$. Then $f_{m,n}=f_{d,n}$.
\end{lem}
\pf
Let $\mathcal A_{m,n}$ denote the family of all $k$-subsets
$\{x_1,\ldots,x_k\}$ of $\mathbb Z_n$ such that
$|x_i-x_j|\notin\{m,2m,\ldots,pm\}$ for all $1\leq i<j\leq k$. Then $f_{m,n}=|\mathcal A_{m,n}|$.
Since $(m,n)=d$, by Lemma \ref{lem:two}, there exist integers $a$ and $b$ such that $(a,n)=1$ and
$am+bn=d$. Let $a^{-1}$ be the inverse of $a\in\mathbb Z_n$.
For any $X=\{x_1,\ldots,x_k\}\in\mathcal A_{m,n}$,
one has $Y=\{ax_1,\ldots,ax_k\}\in A_{d,n}$. Conversely, for any $Y=\{y_1,\ldots,y_k\}\in\mathcal A_{d,n}$,
one can recover $X$ by taking $X=\{a^{-1}y_1,\ldots,a^{-1}y_k\}$.
This proves that $X\mapsto Y$ is a bijection, and therefore $|\mathcal A_{m,n}|=|\mathcal A_{d,n}|$.
\qed

Now we can give a proof of Theorem \ref{thm:msun}.
By Lemma \ref{lem:three}, it suffices to prove it for the case that $n$ is divisible by
$m$.

\medskip
\noindent{\it Proof of Theorem \ref{thm:msun}.} Suppose $n=mn_1$. Let $\mathbb Z_{n,i}=\{i+mj\colon j=0,\ldots,n_1-1\}$.
Then $|\mathbb Z_{n,i}|=n_1$ and $\mathbb Z_n=\biguplus_{i=0}^{m-1}\mathbb Z_{n,i}$.
For any $X=\{x_1,\ldots,x_k\}\subseteq\mathbb Z_n$ and $i=0,\ldots,m-1$, define $X_i=X\cap \mathbb Z_{n,i}$
and $Y_i=\{j\colon j=0,\ldots, n_1-1\text{ and }i+mj\in X_i\}$. Consider $Y_i$ as a subset of $\mathbb Z_{n_1}$.
It is easy to see that $X\in\mathcal A_{m,n}$ if and only if $Y_i\in\mathcal A_{1,n_1}$ for all $i=0,\ldots,m-1$.
Let $|Y_i|=|X_i|=k_i$.  By the aforementioned Kaplansky's result, we have the following expression:
\begin{align}
f_{m,n}=\sum_{k_1+\cdots+k_m=k}\prod_{i=1}^m\frac{n_1}{n_1-pk_i}{n_1-pk_i\choose k_i}.
\label{eq:fmn}
\end{align}
Note that $n\geq mpk+1$, i.e., $n_1\geq pk+1$, the above expression is always well-defined. Finally,
by repeatedly using Rothe's identity
\begin{align*}
\sum_{k=0}^{n}\frac{xy}{(x+kz)(y+(n-k)z)}{x+kz\choose k}{y+(n-k)z\choose n-k}=
\frac{x+y}{x+y+nz}{x+y+nz\choose n}
\end{align*}
(see \cite{Blackwell, Gould, Guo, Rothe}), one sees that
$$f_{m,n}=\frac{n}{n-pk}{n-pk\choose k}.$$
\qed

\noindent{\it Remark.} The idea of writing $\mathbb Z_n$ as a union of some pairwise non-intersecting subsets
is the same as that in \cite[Section 2]{MS}. However, we are unable to obtain such an expression for $f_{m,n}$
if $n\not\equiv 0\pmod m$, as mentioned by Mansour and Sun \cite{MS}.
This is why we need to establish Lemma \ref{lem:three}. Our proof may be deemed as a semi-bijective proof, and
finding a purely bijective proof of Theorem \ref{thm:msun} still remains open.


\end{document}